\documentclass[11pt]{article}

\usepackage{amsmath}
\usepackage{amssymb}
\usepackage{eucal}

\begin{document}

\baselineskip 16pt

\title{Finite groups with Frobenius normalizer condition for non-normal
 primary subgroups  \thanks{Research  is supported by
 China Scholarship Council,  NNSF of
China(11771409) and Wu Wen-Tsun Key Laboratory of Mathematics of Chinese Academy of Sciences.}
}

\author{Zhang Chi , Wenbin Guo\thanks{Corresponding author}\\
{\small Department of Mathematics, University of Science and
Technology of China,}\\ {\small Hefei 230026, P. R. China}\\
{\small E-mail:
zcqxj32@mail.ustc.edu.cn, wbguo@ustc.edu.cn}}

\date{}
\maketitle

\date{}
\maketitle

\begin{abstract}   A finite group  $P$ is said to be
\emph{primary }  if $|P|=p^{a}$ for some prime $p$.
 We say a primary subgroup  $P$ of a finite group $G$
 satisfies   the
\emph{Frobenius normalizer condition} in $G$
         if  $N_{G}(P)/C_{G}(P)$  is a $p$-group provided $P$ is $p$-group.

In this paper, we determine the structure of a finite group $G$  in which
  every non-subnormal primary
 subgroup
 satisfies    the   Frobenius normalized condition.  In particular, we
prove that if    every non-normal primary
 subgroup of $G$
 satisfies    the   Frobenius condition, then $G/F(G)$ is cyclic and every
maximal non-normal nilpotent subgroup $U$ of $G$ with $F(G)U=G$  is a Carter subgroup
of $G$.

\end{abstract}


\let\thefootnote\thefootnoteorig

\section{Introduction}

Throughout this paper, all groups are finite and $G$ always denotes
a finite group. Moreover,  $\mathbb{P}$ is the set of all    primes,
$p\in \pi \subseteq  \Bbb{P}$ and  $\pi' =  \Bbb{P} \setminus \pi$; $p'=\mathbb{P}\setminus \{p\}$. If $n$ is an integer,
 the symbol $\pi (n)$
denotes the
 set of all primes dividing $n$; as usual,  $\pi (G)=\pi (|G|)$, the set of all
 primes dividing the order of $G$. Recall also that $G$ is said to be
\emph{primary }  if $|G|=p^{a}$ for some prime $p$.  A normal subgroup $N$
of $G$ is said to be \emph{hypercentral in $G$} if either $N=1$ or every
chief factor $H/K$ of $G$ below $N$ is central, that is, $C_{G}(H/K)=G$.
The  product of all hypercentral subgroups of $G$  is called the
\emph{hypercentre}  of $G$ and denoted be $Z_{\infty}(G)$.

The nature of the embedding  of primary subgroups in the group has a
significant effect on the structure of this group (see
 the books \cite{prod, GuoII} and the surveys in  \cite{Ore, commun}).
Recall, for example,  that by  the well-known  Frobenius theorem
\cite[IV, Satz 5.8]{hupp},
 $G$ is $p$-nilpotent
 if and only if for every $p$-subgroup $P$ of $G$
the section
 $N_{G}(P)/C_{G}(P)$  is a $p$-group too.

{\bf Definition 1.1.}      We say that a primary subgroup  $P$ of $G$
 satisfies   the
\emph{Frobenius normalizer condition} in $G$
         if  $N_{G}(P)/C_{G}(P)$  is a $p$-group provided $P$ is a $p$-group.

Thus $G$ is nilpotent if and only if every  primary
 subgroup  of $G$
 satisfies   the
Frobenius normalizer condition in $G$.

Before continuing, consider the following example.

{\bf Example 1.2.} (i) Let $p  >  q  > r > t$ be primes, where $q$ divides $p-1$ and
 $t$ divides $r-1$.
Let  $C_{p}\rtimes
C_{q}$  be
  a non-abelian group of order $pq$, $C_{r}\rtimes C_{t}$  be
  a non-abelian group of order $rt$,  and  let $G=(C_{r}\rtimes C_{t}) \times (C_{p}\rtimes
C_{q})$.  Then  every non-normal primary  subgroup of $G$
 satisfies    the
Frobenius normalizer condition in $G$.

(ii) Recall that $G$ is called    \emph{semi-nilpotent} (see
  Section 7  in
 \cite[Ch. 4]{We} if
the normalizer of every non-normal nilpotent subgroup of $G$ is nilpotent.
It is clear that in every semi-nilpotent group $G$, every non-normal primary
 subgroup  satisfies   the Frobenius normalizer condition in $G$. The converse is not true
 in general case (see the group $G$ in Part (i)).

The Example  1.2(i) shows that groups in which every non-normal primary
 subgroup
 satisfies   the
Frobenius normalizer condition may be non-nilpotent.
Nevertheless,  our first  result shows that the groups with such a property
 have the structure very close
to  the structure of nilpotent groups.

{\bf Theorem 1.3. } {\sl If    every non-normal primary
 subgroup of $G$
 satisfies    the   Frobenius normalizer  condition, then  $G/F(G)$ is cyclic and
 all maximal nilpotent
subgroups $U$ of $G$ with $F(G)U=G$ are Carter subgroups of $G$.  }

We prove  Theorem 1.3 being based on the following our general result.

{\bf Theorem 1.4. } {\sl Suppose that   every non-subnormal  primary
 subgroup of $G$
 satisfies    the   Frobenius normalizer condition. Then $G$ is either of the following type:}

(a) {\sl $G$ is nilpotent.}

(b) {\sl $G$ is soluble and $G$ has a Sylow basis $P_{1}, \ldots ,  P_{t}$
 such that:}

(i) {\sl  For some $1\leq r < t$ the subgroups
$P_{1}, \ldots ,  P_{r}$ are normal in $G$,  $P_{i}$ is not   normal in
$G$ for all $i > r$
 and $E= P_{r+1} \cdots  P_{t} $  is nilpotent.}

(ii) {\sl $F(G)$ is a maximal
 nilpotent subgroup of $G$  and
  $F(G)=F_{0}(G)Z_{\infty}(G)$, where
 $F_{0}(G)=P_{1} \cdots   P_{r}$.}

(iii) {\sl $N_{G}(E)$ is a Carter subgroup  of $G$.}

(iv) {\sl $V_{G}= Z_{\infty}(G)$
for every maximal nilpotent subgroup $V$ of $G$ such that  $G=F(G)V$.
}

{\sl Conversely, if $G$ is a group of  type (a) or (b), then
every non-subnormal  primary       subgroup of $G$
 satisfies    the   Frobenius normalizer condition. }

In view of Example  1.2(ii) we get from Theorems  1.4  the following

{\bf Corollary 1.5} (See   Theorem 7.6 in \cite[Ch. 4]{We}).  {\sl If $G$ is
 semi-nilpotent and $F_{0}(G)$ denotes the  product of its normal
 Sylow subgroups, then
$G/F_{0}(G)$ is nilpotent. }

{\bf Corollary 1.6} (See   Theorem 7.8 in \cite[Ch. 4]{We}).  {\sl If $G$ is
semi-nilpotent, then: }

(a) {\sl $F(G)$ is a maximal nilpotent subgroup of $G$.}

(b) {\sl If $U$ is  a maximal nilpotent subgroup of $G$ and $U$ is not normal in $G$, then
 $U_{G}=Z_{\infty}(G)$.}

{\bf Proof.}  If $U$ is  a maximal nilpotent subgroup of $G$ and $U$ is not normal in $G$, then
$N_{G}(U)=U$, so $U$ is a Carter subgroup of $G$. Hence we get both the Statements (a) and (b)
 from Theorem 1.4.

 From Theorems 1.3 we get   the following

{\bf Corollary 1.7} (Chih-Han Sah  \cite{Sah}).  {\sl If $G$ is
semi-nilpotent, then  $G/F(G)$ is cyclic.    }

\section{Proof of Theorem 1.4}

Recall that  $G$ is called a \emph{Schmidt group} if $G$ is not nilpotent
but
 every proper subgroup of $G$ is nilpotent.

 {\bf Lemma 2.1} (See  \cite[Ch. 1, Proposition 1.9]{{GuoII}}).  {\sl If
$G$ is a
 Schmidt group, then $G=P\rtimes Q$, where    $P=G'$
 is a Sylow $p$-subgroup of $G$ and $Q$ is a
 Sylow $q$-subgroup of $G$ for some primes $p\ne q$.    }

$G$ is said to be \emph{$p$-decomposable} if $G=O_{p}(G)\times
O_{p'}(G)$.

 {\bf Lemma  2.2.}  {\sl Suppose that $G$ is soluble and
 let
$P_{1}, \ldots , P_{t}$ be  a Sylow basis of $G$,
 where $P_{i}$ is a
 $p _{i}$-group for all $i$.  If $N_{G}(P_{i})$  is
 a $p _{i}$-decomposable for  all  $i=1, \ldots t$,
 then $G$ is nilpotent.  }

 {\bf Proof. }   Let $R$ be a minimal normal
subgroup of $G$. Then $R$ is a $p _{k}$-group for some $k$ since $G$
is  soluble by hypothesis.
Moreover,  $P_{1}R/R, \ldots , P_{t}R/R$ is  a Sylow basis of $G/R$.
It is clear also that   $N_{G}(P_{i})R=N_{G}(P_{i}R)$, so
 $$N_{G/R}(P_{i}R/R)=N_{G}(P_{i})R/R \simeq N_{G}(P_{i})/N_{G}(P_{i})\cap R$$
is $p_{i}$-decomposable for all $i$. Therefore the hypothesis holds for $G/R$, so $G/R$
is  nilpotent  by induction. Hence $P_{k}/R$ is normal in $G/R$, so
 $G=N_{G}(P_{k})$  is $p _{k}$-decomposable. Let $N$ be a minimal
normal subgroup of $G$ contained in $O_{p_{k}'}(G)$. Then  $G/N$ is
nilpotent, so $G\simeq G/1=G/R\cap N$ is  nilpotent.
  The lemma    is proved.

{\bf Lemma 2.3.} {\sl  Let $H$, $E$ and $N$ be subgroups of $G$, where
 $H$ is a   $p$-subgroup of $E$ and $N$ is a  soluble normal subgroup of  $G$.
 If $N_{G}(H)/C_{G}(H)$ is a  $p$-group   and either $N\leq H$ or $(|N|, |H|)=1$, then
 $N_{G/N}(NH/N)/C_{G/N}(HN/N)$  and $N_{E}(H)/C_{E}(H)$ are $p$-groups.}

{\bf Proof. }  First  assume that $(|N|, |H|)=1$. Then $N_{G}(NH)=NN_{G}(H)$,    so
 $$N_{G/N}(NH/N)/(C_{G}(H)N/N)= (N_{G}(NH)/N)/(C_{G}(H)N/N)=
  (N_{G}(H)N/N)/(C_{G}(H)N/N)$$$$
\simeq N_{G}(H)N/C_{G}(H)N\simeq N_{G}(H)/(N_{G}(H)\cap NC_{G}(H))=
N_{G}(H)/C_{G}(H)(N_{G}(H) \cap N) $$$$\simeq (N_{G}(H)/C_{G}(H))/(C_{G}(H)(N_{G}(H)
 \cap N) /C_{G}(H))$$  is  $p$-group, where
$C_{G}(H)N/N\leq C_{G/N}(NH/N)$. Hence  $N_{G/N}(NH/N)/C_{G/N}(NH/N)$  is a
$p$-group.
  Similarly one can shows that  $N_{G/N}(H/N)/C_{G/N}(H/N)$  is a
$p$-group in the case when $N\leq H$.

Finally,
$$N_{E}(H)/C_{E}(H)=
 (N_{G}(H)\cap E)/(C_{G}(H)\cap E)\simeq (N_{G}(H)\cap
E)C_{G}(H)/C_{G}(H)$$ is a   $p$-group.
The lemma is proved.

 {\bf Lemma 2.4.} {\sl  If $A$ is a subnormal nilpotent $\pi$-subgroup of
$G$, then $A\leq O_{\pi}(G)\cap F(G)$.}

{\bf Proof. }   There is a subgroup chain  $A=A_{0} \leq A_{1} \leq \cdots \leq
A_{n}=G$  such that  $A_{i-1}\trianglelefteq A_{i}$   for all $i=1, \ldots , n$.
We can assume  without loss of generality that $M=A_{n-1}  < G$.  Then by
induction we have that $A\leq O_{\pi}(M)\cap F(M)$. On the other hand, the
subgroups $O_{\pi}(M)$ and $F(M)$ are normal in $G$ since they are
characteristic in $M$, so  $O_{\pi}(M)\leq O_{\pi}(G)$ and $F(M)\leq
F(G)$. The lemma is proved.

{\bf Proposition 2.5.} {\sl  Suppose that for any  Sylow $p$-subgroup
 $P$  of $G$,     every
non-subnormal subgroup  $H$ of $G$ contained in the focal subgroup $G'\cap P$
 satisfies    the   Frobenius  normalizer
condition in $G$.   Then $G$ is $p$-soluble.}

{\bf Proof. } Assume that this proposition is false and let $G$ be
 a counterexample of minimal order. Then $p\in \pi (G)$.

(1) {\sl If $A$ is any $p$-closed Schmidt
subgroup of $G$,  where $p$ divides $|A|$, then
for some minimal normal subgroup $R$ of $G$ we have $R\leq O_{p}(G)$. }

Let $A_{p}$ be the normal Sylow subgroup of $A$. Then $A_{p}=A'\leq G'\cap
G_{p}$ for some Sylow $p$-subgroup $G_{p}$ of $G$ by Lemma 2.1.  On the other hand,
$A_{p}$
 does not satisfy the
 Frobenius normalizer condition in $G$ since $A\leq N_{G}(A_{p}) $   and
$A$ is not $p$-nilpotent. Therefore $A_{p}$ is subnormal in $G$ by
hypothesis, so $A_{p}\leq O_{p}(G)$ by Lemma 2.4.  Hence we have (1).

(2)  {\sl  Every maximal subgroup of $G$ is $p$-soluble.}

First we show that the hypothesis holds for every subgroup $E$ of $G$.
Indeed, let $p\in   \pi(E)$ and let $H$ be any  non-subnormal  $p$-subgroup of $E$ such that
for some Sylow $p$-subgroup $E_{p}$ of $E$ we have $H\leq E'\cap E_{p}$.
Then for some Sylow $p$-subgroup $G_{p}$ of $G$ we have $E_{p}\leq G_{p}$.
On the other hand, $H$ is not subnormal in $G$ by \cite[Ch. A, Lemma 14.1(a)]{DH} since
 it  not subnormal in $E$  and
$H\leq G'\cap G_{p}$ since $E'\leq G'$. Hence  $H$
  satisfies    the   Frobenius normalizer  condition
 in $G$ by hypothesis and hence $H$
  satisfies    the   Frobenius normalizer condition in $E$ by Lemma 2.3.  Therefore
the hypothesis holds for $E$, so the choice of $G$ implies
that every maximal subgroup $M$ of $G$ is $p$-soluble.

(3)  {\sl $F(G)=\Phi (G)$ and $G/F(G)$ is a non-abelian simple group of order divisible by
 $p$. Hence $G'=G$.}

 Let $N$ be any normal proper subgroup of $G$. Then $N$ is $p$-soluble by Claim (2).
Suppose that $N\nleq M$ for some maximal subgroup $M$ of $G$. Then
$G/N\simeq M/M\cap N$ is $p$-soluble by Claim (2) and so $G$ is $p$-soluble,
 contrary to our hypothesis about $G$.
Therefore $N\leq \Phi (G)\leq F(G)$, so $F(G)=\Phi (G)$ and $G/F(G)$ is a
non-abelian simple group of order divisible by $p$.   Hence we have (3).

{\sl The final contradiction.}  Since $G$ is not $p$-soluble, it is not
$p$-nilpotent. Hence $G$ has a subgroup $A$ such that $A$ is not
$p$-nilpotent but every maximal subgroup of $A$ is
$p$-nilpotent. Then $A$ is a $p$-closed Schmidt group with $p\in \pi (A)$
by \cite[Ch. IV,
Satz 5.4]{hupp}.  Then
for some minimal normal subgroup $N$ of $G$ we have $N\leq O_{p}(G)$ by
Claim (1).

 Let  $H/N$ be any  non-subnormal  $p$-subgroup of $G/N$ such that
for some Sylow $p$-subgroup $P/N$ of $G/N$ we have $H/N\leq  (G/N)'\cap
(P/N)$. Then $P$ is a Sylow $p$-subgroup of $G$ and, by Claim (3),  $H$ is a
non-subnormal  $p$-subgroup of $G$ such that  $H\leq P=P\cap G'$.
 Therefore   $H$
  satisfies    the   Frobenius normalizer condition in $G$ by hypothesis and hence
 $H/N$
  satisfies    the   Frobenius normalizer condition in $G/N$ by Lemma 2.3.  Therefore
the hypothesis holds for $G/N$. The choice of $G$ implies
that $G/N$ is $p$-soluble and so  $G$ is $p$-soluble.
The proposition is proved.

We use  in our proof the following properties of the subgroup $Z_{\infty}(G)$.

{\bf Lemma 2.6 }  (See Theorem 2.6 in \cite[Ch. 1]{GuoII}).   {\sl  Let   $Z=Z_{\infty}(G)$.
 Let  $A$, $B$   and $N$ be subgroups of $G$, where
$N$ is normal in $G$.}

(1) {\sl If $ N\leq Z$, then $Z/N= Z_{\infty}(G/N)$.}

(2)   {\sl If $A$ is nilpotent, then $ZA$   is also nilpotent. Hence
$Z$ is contained in each maximal nilpotent subgroup of $G$.}

(3) {\sl If $G/Z$ is nilpotent, then $G$   is also nilpotent.}

  The following Lemma is a corollary of Lemma 2.6(3) and \cite[Ch. A, Theorem 9.3(c)]{DH}.

 {\bf Lemma 2.7.} {\sl  $F(G)/\Phi (G)=F(G/\Phi(G))$ and
$F(G)/Z_{\infty}(G)=F(G/Z_{\infty}(G))$.}

 Let $\phi$ be some linear ordering on  $\mathbb P$.
The record $p\phi q$  means that $p$  precedes $q$ in $\phi$ and $p\ne q$.
Recall that
a group $G$ of order
$p_1^{\alpha _1}p_2^{\alpha _2}\ldots
p_n^{\alpha _n}$
is called \emph{ $\phi$-dispersed} (Baer \cite{baer1}) whenever
$p_{1}\phi p_{2}\phi \ldots \phi p_{n}$
and for every $i$ there is
a normal subgroup of $G$ of order
$p_1^{\alpha _1}p_2^{\alpha _2}\ldots
p_i^{\alpha _i}$.

{\bf Proposition 2.8.} {\sl Suppose that $G$ is soluble.}

(a) {\sl   If every non-normal Sylow subgroup $P$ of $G$ with
 $P\cap G'\ne 1$  satisfies    the   Frobenius normalizer
condition in $G$, then $G$  $\phi$-dispersed for some   linear ordering
$\phi$ on  $\mathbb P$.}

(b)  {\sl If $G$ is not nilpotent but every non-normal Sylow subgroup  of $G$
  satisfies    the   Frobenius
normalizer condition in $G$, then the following conditions holds:}

(i) {\sl  $G$  has a Sylow basis $P_{1}, \ldots ,  P_{t}$
 such that for some $1\leq r < t$ the subgroups
$P_{1}, \ldots ,  P_{r}$ are normal in $G$,  $P_{i}$ is not   normal in
$G$ for all $i > r$
 and $E= P_{r+1} \cdots  P_{t} $  is nilpotent.}

(ii) {\sl $F(G)$ is a maximal
 nilpotent subgroup of $G$  and
  $F(G)=F_{0}(G)Z_{\infty}(G)$, where
 $F_{0}(G)=P_{1} \cdots   P_{r}$.}

(iii) {\sl $N_{G}(E)$ is a Carter subgroup  of $G$.}

(iv) {\sl $V_{G}= Z_{\sigma}(G)$
for every maximal nilpotent subgroup $V$ of $G$ such that  $G=F(G)V$.
}

  {\bf Proof. }  Since $G$ is soluble, it has a Sylow basis $P_{1},
\ldots, P_{t}$.
We can assume without loss of generality that $P_{i}$ is a
$p_{i}$-group for all $i=1, \ldots t.$

 (b)   Assume that this assertion  is false and let $G$ be
 a counterexample of minimal order.

 (1) {\sl  If  $N$ is a  minimal normal subgroup of $G$, then
the  conclusions of the proposition  hold  for
$G/N$.  }

Since $G$ is soluble, $N$ is  primary.
Let $p\in \pi (G/N)$. Suppose that  a Sylow $p$-subgroup $P/N$  of   $G/N$
is not normal  $G/N$. Then $P/N=G_{p}N/N$ for some Sylow $p$-subgroup
$G_{p}$ of $G$. Moreover, \cite[Ch. A, Lemma 14.1(b)]{DH} implies that $G_{p}$
 is not normal in $G$. Hence  $G_{p}$
  satisfies    the   Frobenius condition in $G$ by hypothesis and hence $P/N=G_{p}N/N$
  satisfies    the   Frobenius condition in $G/N$ by Lemma 2.3.
 Therefore the hypothesis holds for $G/N$, so we
 have (1) by the choice of $G$.

(2) {\sl If $P_{i}$ is not normal  in $G$, then $N_{G}(P_{i})$ is $p _{i}$-decomposable. }

Since
$N_{G}(P_{i})/C_{G}(P_{i})$ is a  $p_{i}$-group by hypothesis,  a Hall
 $p_{i}'$-subgroup $V$ of $N_{G}(P_{i})$ is normal in $N_{G}(P_{i})$.
Hence we have (2).

(3) {\sl Statement (b)(i) holds for $G$.}

Since $G$ is not nilpotent, Lemma 2.2 and Claim (2) imply that
   for some $1\leq r < t$ the subgroups
$P_{1}, \ldots ,  P_{r}$ are normal in $G$ and   $P_{i}$ is not   normal in
$G$ for all $i > r$. Then $E= P_{r+1} \cdots  P_{t} $  is nilpotent by Claim (2)
 and Lemma 2.3.  Hence we have  (3).

(4) {\sl Every subgroup $V$
 of $G$ containing $ F(G)$ is subnormal in $G$, so
 $F(V)=F(G)$}.

First note that $G/F(G)=EF(G)/F(G)\simeq E/E\cap F(G)$ is nilpotent by Claim (3).
Therefore $V/F(G)$ is subnormal in $G/F(G)$, which implies that $V$ is subnormal in $G$.
Then  $F(V)\leq F(G)$ by Lemma 2.5. On the other hand, $F(G)\leq F(V)$ and
so   $F(V)=F(G)$.

(5)  {\sl Statement (b)(ii)  holds for $G$.}

First note that  every nilpotent subgroup $V$ of $G$    containing $F(G)$
is subnormal in $G$ by Claim (4), so $V\leq F(G)$ by Lemma 2.5.
Therefore $F(G)$  is a maximal
 nilpotent subgroup of $G$.

In fact, $F(G)=F_{0}(G)\times O_{p_{i_{1}}}(G)\times \cdots \times
 O_{p_{i_{m}}} (G)$ for some $i_{1}, \ldots , i_{m} \subseteq \{r+1,
\ldots , t\}$.
Moreover, in view of Claim (6), we get   that $G/C_{G}(O_{p _{i_{k}}}(G))$  is a
$p_{i_{k}}$-group    for every  $k=1,  \ldots m$ and so  for every chief
factor  $H/K$ of $G$ with   $H\leq O_{p _{i_{k}}}(G)$ we have
$C_{G}(H/K)=G$ since $O_{p _{i_{k}}}(G/C_{G}(H/K))=1$ by \cite[Appendixes, Corollary 6.4]{We}.
 Therefore   $O_{p _{i_{k}}}(G)\leq Z_{\infty}(G)$.
  Hence
  $F(G)=F_{0 }(G)Z_{\infty}(G)$.

(6)  {\sl Statement (b)(iii)  holds for $G$.}

The subgroup $N_{G}(P_{i})$ is $p_{i}$-decomposable for all $i > r$ by Claim
(2). Therefore, by Claim (3),  $$N_{G}(E)=
 N_{G}(P_{r+1}\times \cdots \times P_{t})=N_{G}(P_{r+1})\cap \cdots \cap
N_{G}(P_{t})$$ is nilpotent. On the other hand, $N_{G}(N_{G}(E))=N_{G}(E)$  since $G$ is soluble.
 Hence we have (6).

(7)   Statement (b)(iv)  holds for $G$

First we show that $U_{G}\leq  Z_{\infty}(G)$ for every nilpotent subgroup $U$
 of $G$               such that  $G=F(G)U$.  Suppose that this is false. Then $U_{G}\ne 1$.
 Let $R$ be a minimal
normal subgroup of $G$ contained in $U$ and $C=C_{G}(R)$.
  Then $$G/R=(F(G)R/R)(U/R)=F(G/R)(U/R),$$ so $$U_{G}/R=(U/R)_{G/R}\leq
 Z_{\infty}(G/R)$$ by Claim (1). Since $G$ is soluble, $R$
is a $p$-group for some prime $p$.  Moreover, from $G=F(G)U$
we get that for some Hall $p'$-subgroups $E$, $V$ and
$W$ of $G$, of $F(G)$ and of  $U$, respectively, we have $E=VW$.
But $R\leq F(G)\cap U$, where $F(G)$ and $U$ are
nilpotent.   Therefore $E\leq C$, so $R/1$ is central in
$G$. Hence $R\leq Z_{\infty}(G)$ and so $Z_{\infty}(G/R)=Z_{\infty}(G)/R$ by Lemma 2.6(1).
But then $U_{G}\leq Z_{\infty}(G)$. Finally,  $Z_{\infty}(G)\leq  U$ by
Lemma 2.6(2) and so $U_{G}= Z_{\infty}(G)$.

From Claims (3), (5), (6) and  (7)  it follow that all  conclusions of the
proposition  hold for $G$, contrary to the choice of $G$. This final contradiction completes
 the proof of Assertion (b).

  (a) Assume that this assertion  is false and let $G$ be
 a counterexample of minimal order.

First assume that  that for some $p\in
\pi (G)$  and for a  Sylow $p$-subgroup  $G_{p}$ of we have $G_{p}\cap
G'=1$. Then $G$ is  $p$-nilpotent, so $G$ has a  normal $p$-complement
 $E$. Now let $P$ be a  non-normal Sylow subgroup  of $E$ with
 $P\cap E'\ne 1$. Then $P$ is  a  non-normal Sylow subgroup  of $G$ and
$P\cap G'\ne 1$. Hence $P$ satisfies    the   Frobenius normalizer
condition in $G$ and so   $P$ satisfies    the   Frobenius normalizer
condition in $E$ by Lemma 2.3. Therefore the hypothesis holds for $E$, so
$E$ is $\phi$-dispersed for some   linear ordering
$\phi$
on  $\mathbb P$, so $G$ is $\phi _{0}$-dispersed for some   linear ordering
$\phi _{0}$ on  $\mathbb P$, contrary to our assumption about $G$.
Therefore for every Sylow subgroup $G_{p}$ of $G$ we have $G_{p}\cap G'\ne
1$. But in this case Statement (a) is a corollary of Statement (b).

 The proposition  is proved.

{\bf Proof of Theorem 1.4.}    First note   that  if  every non-subnormal  primary
 subgroup of $G$
 satisfies    the   Frobenius normalizer condition, then in view of Propositions 2.5 and 2.8,
 $G$ is one of the types (a) or (b).

To complete the proof of the theorem it is enough to
  show that if  $G$  is a group of type (b), then  every non-subnormal  primary
 subgroup $H$ of $G$
 satisfies    the   Frobenius condition. Let $H$ be  a $p$-group. Then for
some $i > r$  and $x\in G$ we have $H\leq P_{i}^{x}\leq E^{x}$. Let $N=N_{G}(H)$.
 Then $[N\cap
F_{0}(G), H]=1=[V, H]$, where $V$ is a $p$-complement of  $E^{x}$.
Moreover, $G=F_{0}(G)\rtimes E^{x}$ and so
$N=(F_{0}(G)\cap N)(N\cap E^{x})$, which implies that $N_{G}(H)/C_{G}(H)$ is a
$p$-group.
The theorem is proved.

\section{Proof of Theorem 1.3}

                 \

 Suppose that this theorem  is false and let
$G$ be a counterexample of minimal order.  Then $G$ is not nilpotent.
Moreover, Theorem 1.4 implies that $G$ is soluble and it has a Sylow  basis $P_{1}, \ldots ,  P_{t}$
 such
that  for some $1\leq r < t$ the subgroups
$P_{1}, \ldots ,  P_{r}$ are normal in $G$,  $P_{i}$ is not   normal in
$G$ for all $i > r$
 and $E= P_{r+1} \cdots  P_{t} $  is nilpotent.  Let $F_{0}(G)= P_{1}
\cdots   P_{r}$.

(1) {\sl If $A=A_{p}\rtimes A_{q}$ is any $p$-closed Schmidt
subgroup of $G$,  where $p$ divides $|A|$, then $A_{p}$ is normal in $G$.}

The subgroup
$A_{p}$
 does not satisfy the
 Frobenius normalizer condition in $G$ since $A\leq N_{G}(A_{p}) $   and
$A$ is not $p$-nilpotent. Therefore we have (1) by
hypothesis.

(2) {\sl The  conclusions of the theorem  hold for every proper subgroup  of $G$ and
 for every
quotient $G/N$ of $G$, where $N$ is a minimal normal subgroup of $G$}
 (See the proof of Proposition 2.8 and Theorem 1.4).

(3)  {\sl $G/F(G)$ is abelian. }

It is enough to show that $G'$ is  nilpotent. Suppose
that this is false.  Let $R$ be  a minimal normal subgroup of $G$.

(a) {\sl $R=C_{G}(R)=O_{p}(G)=F(G)\nleq \Phi (G)$ for some prime $p$ and $|R|  > p$}.

 From Claim (2) it follows that for every minimal normal subgroup $N$ of $G$,
 $(G/N)'=G'N/N\simeq G'/G'\cap N$
 is nilpotent.

 If $R\ne N$, it follows that  $G'/((G'\cap N)\cap
(G'\cap R))=G'/1$ is  nilpotent. Therefore $R$ is the  unique
minimal normal subgroup of $G$, $R\leq G'$ and, by Lemma 2.7,  $R\nleq \Phi (G)$. Hence
$R=C_{G}(R)=O_{p}(G)=F(G)$   by    Theorem 15.6 in \cite[Ch. A]{DH}, so
  $|R| > p$ since otherwise
       $G/R=G/C_{G}(R)$ is cyclic, which implies that   $G'=R$ is  nilpotent.

(b)  {\sl $G=P_{1} \rtimes P_{2}$, where $R\leq P_{1}=F(G)$ and $P_{2}$ is
 a minimal non-abelian group.}

From      Claim    (a) it follows that $r=1$  and $R\leq P_{1}=F(G)$.

 Now let
$W=P_{1}V$, where $V$ is a  maximal subgroup of $E$. Then $W$ is
subnormal in $G$ and so $F(W)=F(G)=P_{1}$ (see Claim (7) in the proof of
Theorem 1.4).
But then  $W/P_{1}=P_{1}V/P_{1}\simeq V$ is
abelian by Claim (2). Therefore $E$  is not abelian but
 every proper subgroup of $E$  is abelian, so  $E=P_{2}$ since $E$ is
 nilpotent.  Hence we have (b).

(c)  {\sl $P_{1}=R$  is a Sylow $p$-subgroup of $G$ and every  subgroup $H\ne 1$ of
 $P_{2}$ acts irreducibly  on $R$. Hence every proper  subgroup $H$  of $P_{2}$ is cyclic.}

Since  $P_{1}$  is a normal Sylow $p$-subgroup of $G$,
$P_{1}\leq F(G)\leq C_{G}(R)=R$ by Claim (a)   and  \cite[Ch. A, 13.8(b)]{DH} and  so $P_{1}=R$.

Now let $S=RH$.  By
the Maschke's theorem, $R=R_{1}\times \cdots \times R_{s}$, where $R_{i}$
is a minimal normal subgroup of $S$ for all $s$. Then $R=C_{S}(R)=C_{S}(R_{1}) \cap
\cdots \cap C_{S}(R_{n})$. Hence for some $i$
 the subgroup $R_{i}H$ is not nilpotent and so it has
 a Schmidt
subgroup  $A$ such that $1  < A'$ is normal in $G$ by Claim (1). But then
 $R\leq A$. Therefore $i=1$, so we have (c) since $H$ is abelian by Claim (b).

 {\sl The final contradiction for (3). }   Since every maximal subgroup of  $P_{2}$ is cyclic by Claim (c),
 $q=2$ by
 \cite[Ch. 5,  4.3, 4.4]{Gor}. Therefore $|R|=p$, contrary to
Claim (a).  Hence we have (3).

(4)  {\sl   $G/F(G)$ is cyclic.   }

 Assume that this is false.
First we show that  $Z_{\infty}(G)=1=\Phi (F(G))$.
Assume that for some minimal normal subgroup $R$ of $G$ we have either
$R\leq Z_{\infty}(G)$ or $R\leq  \Phi (F(G))$.
Then $F(G/R)=F(G)/R$ by Lemma 2.7, so Claim (2) implies that
 $(G/R)/F(G/R) =(G/R)/(F(G)/R)  \simeq G/F(G)$ is cyclic.
This contradiction shows that  $Z_{\infty}(G)=1=\Phi (F(G))$.   Therefore
$F(G)=F_{0}(G)$ by Theorem 1.4 and also $F_{0}(G)=R_{1} \times \cdots \times
R_{k}$ for some minimal normal subgroups  $R_{1}, \ldots ,
R_{k}$ of $G$    by \cite[Ch.A, Theorem 10.6(c)]{DH}.
 Since $E\simeq G/F(G)$ is abelian by Claim  (3), and
  $G$ is not nilpotent,  there
 is an index $i$ such that $V=R_{i}\rtimes E$ is not
nilpotent. Then $C_{R_{i}}(E)\ne R_{i}$.   By the Maschke's theorem, $R_{i}=
L_{1}\times \cdots \times L_{m}$
  for some minimal normal subgroups $L_{1}, \ldots ,  L_{m}$ of $V$. Then, since
$C_{R_{i}}(E)\ne R_{i}$,  for some $j$ we have $L_{j}\rtimes E\ne L_{j}\times E$. Hence
$L_{j}E$ contains a Schmidt subgroup $A_{p}\rtimes A_{q}$ such that $A_{p}=R_{i}$ by Claim (1), so $m=1$.
But then  $E$ acts irreducible
on $R_{i}$ and  hence $G/F(G)\simeq E$ is cyclic.  This
contradiction completes the proof of the fact that  $G/F(G)$ is cyclic.

(5) {\sl All maximal nilpotent
subgroups $U$ of $G$ with $F(G)U=G$ are Carter subgroups of $G$.}

 Suppose that this is false.  Assume that for some minimal normal
 subgroup $R$ of $G$ we have
$R\leq Z_{\infty}(G)$.
By Lemma 2.6(2),
  $R\leq  U$. On the other,
$U/R$ is a
maximal nilpotent non-normal subgroup of $G/R$  by
Lemma 2.6(3). Hence     Claim (2) implies that
$U/R$ is a  Carter subgroup  $G/R$,  so
 $U$ is a  Carter subgroup  of $G$.  Hence   $Z_{\infty}(G)=1$, so
    Theorem 1.4 implies that $F(G)=
F_{0}(G)=P_{1}\cdots P_{r}$. Hence $E\simeq  G/F_{0}(G)$ is abelian by Claim (3).
Since   $G=F(G)U $,  for some $x$ we have $E^{x}\leq U$. Hence
$U\leq N_{G}(E^{x})$ since $E^{x}$ is a Hall  subgroup of $G$ and so $U=N_{G}(E^{x})$
  is a
 Carter subgroup of $G$ by    Theorem 1.4(b)(iii).
 This contradiction completes
 the proof of  (5).

Claims (4) and (5) show that the conclusions of the theorem hold for $G$,
 contrary to the choice of $G$.    The theorem is proved.

\end{document}